\documentclass[11pt]{article}
\usepackage{amssymb,amsmath}
\usepackage{graphicx}
\usepackage{tabbert}

\newcommand{\Sph}{\mathbb S} 
\newcommand{\hb}{\hat{b}}
\newcommand{\rot}{{\rm rot}}
\newcommand{\tvp}{\tilde{\varphi}}
\newcommand{\tal}{\tilde{\alpha}}

\allowdisplaybreaks[2]

\title{The continuous theory of dislocations for a material containing
  dislocations to one Burgers vector only}

\author{
Hans-Dieter Alber\addtocounter{footnote}{1}\footnote{alber@mathematik.tu-darmstadt.de}
\\
{\small Fachbereich Mathematik, Technische Universit\"at Darmstadt} \\ 
{\small Schlossgartenstr. 7, 64289 Darmstadt, Germany}
}

\date{}
\begin{document}
\maketitle
\begin{abstract}
  We review the continuous theory of dislocations from a mathematical point of
  view using mathematical tools, which were only partly available when the
  theory was developed several decades ago. We define a space of dislocation
  measures, which includes Hausdorff measures representing the dislocation
  measures of single dislocation curves. The evolution equation for
  dislocation measures is defined on this space. It is derived from four basic
  conditions, which must be satisfied by the model.
\end{abstract}
\section{Introduction}\label{S1}

Plastic deformation of metallic bodies is caused by the creation and movement
of dislocations in the crystal lattice of the metallic material. Therefore the
plastic deformation depends on the number of dislocations and on the
restrictions of the dislocation movement by the crystal structure and by the
geometry of the body. Standard phenomenological models for viscoplastic
material behavior do not reflect these material properties depending on the
microstructure of the material. A modelling approach taking this
microstructure into account is the continuous theory of dislocations, which
was developed several decades ago and which is well understood in continuum
mechanics. Under many articles in this field we only mention the
classical expositions \cite{Kroe58,Mura63,Fox66,Kos79} and the articles
\cite{ABR08,Acharya11} containing new developments. This theory has
interesting and difficult mathematical aspects. We hope to make the
mathematical aspects better accessible by reviewing the theory from a
mathematical point of view using mathematical tools, which were only partly
available when the theory was developed. We are convinced that the mathematical
aspects deserve much more investigation and that the theory can be advanced
and simulations based on the theory can be improved by such investigations.

We begin by stating the standard model for the deformation of a viscoplastic
body consisting of a metallic material with dislocations moving only in one
slip plane of the crystal structure. It consists of the equations
\begin{eqnarray}
- \div_x \;T(x,t) &=&0, \label{E1.1}\\
T(x,t) &=& D \Big(\ve \big(\na_x u(x,t)\big)-m \ve_p(x,t)\Big),
  \label{E1.2}\\  
\pa_t \ve_p &=& f \big(m:T(x,t)\big), \label{E1.3}
\end{eqnarray}
which must hold for time $t \geq 0$ and for $x$ varying in the open 
set $\Om \subseteq \R^3$ representing the material points of the body.  The
unknowns are the displacement $u(x,t) \in \R^3$ of the material point $x$ at
time $t$, the Cauchy stress tensor $T(x,t) \in \ES^3$, where $\ES^3$ denotes
the set of symmetric $3 \times 3$-matrices, and the plastic strain $\ve_p(x,t)
\in \R$ along the slip plane. We use the standard notation
\[
\div_x T=\big(\sum\limits_{j=1}^3 \pa_{xj} T_{ij} \big)_{i=1,2,3}.
\]
$\na_x u$ denotes the $3 \times 3$-matrix of first order partial derivations
of $u$, and 
\[
\ve(\na_x u) = \frac{1}{2} \big(\na_x u+(\na_x u)^T \big) \in \ES^3 
\]
is the linear strain tensor. We write $A^T$ for the transpose of a
matrix $A$. The elasticity tensor $D:\ES^3 \to \ES^3$ is a linear,
symmetric, positive definite mapping and the constant matrix $m$ is
given by
\begin{equation}\label{E1.4}
m=\ve (\hat{b} \otimes g)= \frac{1}{2}(\hat{b} \otimes g + g \otimes \hat{b})
 \in \ES^3, 
\end{equation}
where $g \in \R^3$ is the unit vector normal to the slip plane, and
$\hat{b} \in \R^3$ is a unit vector in the direction of plastic slip;
it is therefore a vector in the slip plane. For vectors $a,b \in \R^3$
we write $a \otimes b$ to denote the matrix $(a_i b_j)_{i,j=1,2,3}$.
The scalar product of two $3 \times 3$-matrices $A$, $B$ is denoted by
$A:B= \sum_{i,j=1}^3 a_{ij} b_{ij}$. Finally, $f: \mathcal{D}(f)
\subseteq \R \to \R$ is a given function satisfying 
\[
s \cdot f(s) \geq 0,
\]
for all $s \in \R$. We call $f$ constitutive function. A typical
choice for $f$ is $f(s)=C |s|^{\gamma -1} s$, with constants $C >0$
and $\gamma > 1$.  It is well known that if boundary conditions for
$u$ or $T$ are imposed and if $f$ is a maximal monotone function with
$0 \in \mathcal{D}(f)$ and $f(0)=0$ satisfying suitable growth
conditions, then the initial-boundary value problem to \eq{1.1} --
\eq{1.3} has a unique solution. This is proved for example in
\cite{AlChel07}. 

The constitutive equation \eq{1.3} is an ordinary differential
equation in time and does not reflect material behavior caused by the
dislocation microstructure.  The field equations of the continuous
theory of dislocations do reflect this microstructure. To derive these
field equations we start in Section~\ref{S2} by discussing the
Volterra model for dislocation curves in linear elasticity and by
defining the space of dislocation measures. The evolution equation for
dislocation measures is derived in Section~\ref{S3} from four basic
principles. In Sections~\ref{S3.3} we discuss the restrictions on the
model equations following from the incompressibility constraint for
the plastic part of the strain tensor and we obtain the final field
equations. The simplification following from the assumption that
dislocations move in slip planes is discussed in Section~\ref{S4}. At
the end of that section we compare the field equations thus derived
to the standard equations \eq{1.1} -- \eq{1.3}.

As usual in the derivation of model equations, we cannot define the
function spaces, to which the solutions of the model equations belong,
with the same level of rigour as in investigations of existence. This
concerns in particular the space of dislocation measures. We must
assume that the solutions have certain properties; the exact
properties can be determined only after the model equations are known.

Several  technical proofs are not included in this article to keep the
length acceptable. These proofs will be published elsewhere.
\section{The boundary value problem for the stress field of a
  dislocation curve}\label{S2} 
In this section we review the Volterra model for dislocation curves
within the linear theory of elasticity, cf. \cite{Kos79}. This model
suggests the definition of the space $\EM_d(\Om)$ of dislocation
measures, which is given at the end of the section. 

Let $\ell$ be a closed curve in the elastic body $\ov{\Om}$, which
represents a dislocation curve. It is allowed that a part of the curve
belongs to the boundary $\pa\Om$. We assume that an arc length
parametrization $s \mapsto y(s)$ is given, which is continuously
differentiable with the exception of at most finitely many points.
This parametrization defines a unit tangent vector field $\tau=
\frac{\rm d}{{\rm d} s} y$ along $\ell$.  Let $\Sigma$ be a surface in
$\Om$ with $\pa \Sigma=\ell$ and let $n:\Sigma \to \R^3$ be a continuous
normal vector field on $\Sigma$. We choose $n$ such that at $x \in
\pa\Sigma$ the vector $n(x) \times \tau(x)$ points into the surface
$\Sigma$. For functions $v$ defined on $\Om \setminus (\Sigma \cup
\ell)$ and for $x \in \Sigma$ we use the notations
\begin{equation}\label{E2.1}
v^\pm (x)= \lim\limits_{\eta \searrow 0} v \big(x \pm \eta\, n(x)\big), \quad 
  [v]_{\Sigma} (x)=v^+(x)-v^- (x). 
\end{equation}
With the dislocation curve $\ell$ we associate a fixed vector $b \in
\R^3$, $b \not= 0$, the Burgers vector of the dislocation curve. We use
the notation  
\[
\hat{b}= \frac{b}{|b|}
\]
for the unit vector in direction of $b$.

To compute the stress field generated by the dislocation curve in the body
$\Om$, consider the boundary and transmission problem for the displacement
field $u: \ov{\Om} \setminus (\Sigma \cup \ell) \to \R^3$ and the Cauchy
stress tensor field $T: \ov{\Om} \setminus (\Sigma \cup \ell) \to \ES^3$:
\begin{eqnarray}
-\div\, T &=& 0, \label{E2.2}\\
T &=& D \ve (\na u), \label{E2.3}\\
{[u]}_\Sigma &=& -b, \label{E2.4}\\
{[T]}_\Sigma n &=& 0, \label{E2.5}\\
T \rain{{\pa\Om}} n_B &=& \gamma. \label{E2.6}
\end{eqnarray}
The elasticity equations \eq{2.2} and \eq{2.3} must hold on $\Om
\setminus (\Sigma \cup \ell)$, equations \eq{2.4} and \eq{2.5} are
jump relations on $\Sigma$, and \eq{2.6} is the boundary condition on
$\pa\Om$, where $n_B(x) \in \R^3$ denotes the unit normal to $\pa\Om$
at $x \in \pa\Om$ pointing to the exterior of $\Om$ and $\gamma:
\pa\Om \to \R^3$ are the given boundary data. This problem describes
the displacement and stress fields in an elastic body, which is cut
along the surface $\Sigma$. After cutting, the two boundary parts
created by the cutting are displaced elastically against one another
by the vector $b$ and glued together again. Since the length of $b$ is
approximately equal to the lattice constant of the crystal lattice,
after this procedure the atoms on both sides of $\Sigma$ are again in
the right positions to form an elastically stressed, but otherwise
perfect crystal at $\Sigma$. The crystal is disturbed only along the
boundary $\ell$ of $\Sigma$, making \eq{2.2} -- \eq{2.6} a model for
the dislocation curve $\ell$. The stress field generated by this
dislocation curve is given by $T$. This stress field will have a
singularity along $\ell$.

Our goal is to start from this model and to generalize it to a model
for bodies containing an array of dislocations described by a
dislocation density. To do this rigorously, it would be necessary to
solve the problem \eq{2.2} -- \eq{2.6} in a suitable function space and
to study the singularity of $T$ along $\ell$. This difficult task is
out of the scope of this article. Instead, we consider the dislocation
problem in a different situation with a simpler geometry, where an
explicit solution is known.

Namely, we assume that $\Om$ is equal to $\R^3$, that the dislocation
curve $\ell$ is equal to the $x_3$--axis, and that the material is
isotropically elastic, which means that the elasticity tensor is given
by
\begin{equation}\label{E2.7}
D_I \ve = \la\, {\rm trace}(\ve)I + 2 \mu \ve,
\end{equation}
for all $\ve \in \ES^3$, where $I$ is the identity matrix and where
$\la$, $\mu$ are material constants satisfying $\mu > 0$ and $3\la + 2
\mu > 0$.  

To formulate the boundary value problem for the stress field in this
situation we  assume that the Burgers vector is of the form $b =
(b_1,0,b_3)$. This can always be achieved by rotation of the
coordinate system around the $x_3$--axis. We define the tangential
vector field $\tau$ along the line $\ell$ by $\tau(x) = e_3 =
(0,0,1)$. We are free to choose for $\Si$ any half plane with boundary
$\ell$. Therefore we take     
\[
\Sigma=\{(x_1,0, x_3) \; \mid \; x_1>0, \; x_3 \in \R\}
\]
and define the unit normal vector field $n$ on $\Sigma$ by
$n(x)=e_2$. With this definition we obtain for $x \in \ell$ that the
vector $n(x) \ti \tau(x) = e_2 \ti e_3 = e_1$ points into $\Si$, hence
our requirement for the orientation of $n$ is satisfied. Moreover,
from \eq{2.1} we obtain for functions $v$ defined on $\R^3 \setminus
(\Sigma \cup \ell)$ and $x \in \Sigma$ that 
\[
 v^\pm (x)= \lim_{\eta \searrow 0} v (x_1, \pm \eta, x_3),\quad
 [v]_\Sigma \, (x)=\lim_{\eta \searrow 0} \big(v(x_1,\eta,x_3) -
 v(x_1,-\eta,x_3)\big).  
\]
The problem corresponding to \eq{2.2} -- \eq{2.6} for the stress field
generated by the dislocation line $\ell$ consists of the equations 
\begin{eqnarray}
-\div\, T(x) &=& 0, \label{E2.8} \\
T(x) &=& D_I \ve \big(\na u (x)\big), \label{E2.9} \\ 
{[u]}_{\Sigma} (x) &=& -(b_1,0,b_3), \label{E2.10}\\
\big([T]_\Sigma(x)\big) e_2 &=& 0, \label{E2.11}
\end{eqnarray}
where the first two equations must hold for $x \in \R^3 \setminus (\Si
\cup \ell)$ and the last two for $x \in \Si$. A solution of this
problem is given in \cite[pp. 49 -- 53]{Kos79}. In the following theorem we
state this solution and give some additional properties of it. We use the
notation $x = (x',x_3) \in \R^3$ with $x' = (x_1,x_2)$. 
\begin{theo}\label{T2.1}
Let $u = (u_1,u_2,u_3):\R^3 \setminus \ov{\Si} \ra \R^3$ be defined by 
\begin{eqnarray}
\label{E2.12} u_1(x) &= & \frac{b_1}{2\pi} \Ba \arctan \frac{x_2}{x_1} +
\frac{1}{2(1-\nu)}\: \frac{x_1x_2}{r^2}\Bz, \\
\label{E2.13} u_2(x) & = & - \frac{b_1}{4\pi (1-\nu)} \, \Ba (1-2\nu)\log r +
\frac{x_1^2}{r^2}\Bz, \\
\label{E2.14} u_3(x) & = & \frac{b_3}{2\pi} \, \arctan \frac{x_2}{x_1}\,,
\end{eqnarray}
where $\nu = \frac{\la}{2(\la + \mu)} $ is Poisson's ratio
and $r^2 = |x'|^2 = x_1^2 + x_2^2$. Also, let the symmetric tensor
function $T = (T_{ij})_{i,j=1,\ldots, 3}: \R^3 \setminus \ov{\Si} \ra
\ES^3$ be given by  
\begin{eqnarray}
 \label{E2.15}&&\hspace*{-2cm} T_{11} = -D_1\, \frac{x_2(3x_1^2 +
 x_2^2)}{r^4}\,, \; T_{12} = D_1 \, 
\frac{x_1(x_1^2 - x_2^2)}{r^4}\,, \;  T_{13} = -D_2\, \frac{x_2}{r^2}\,
,\\
\label{E2.16}&& \hspace*{-2cm}T_{22} =
D_1\, \frac{x_2(x_1^2-x_2^2)}{r^4}\,,\hph{-3} \;  T_{23} = D_2\, 
\frac{x_1}{r^2}\,,\hph{\frac{x_4(-x_2^2)}{}} \;  T_{33} = -2\nu
D_1\, \frac{x_2}{r^2}\,, 
\end{eqnarray}
with $D_1 = \frac{\mu b_1}{2\pi(1-\nu)}\,,$ $D_2 = \frac{\mu b_3}{2\pi}\,$. 
Then the function $(u,T)$ solves \eq{2.8} -- \eq{2.11}, has the asymptotic
behavior 
\begin{equation}\label{E2.17} 
\na u(x) = O \left(\frac{1}{|x'|} \right),\quad T(x) = O\left
( \frac{1}{|x'|}\right ), \quad \mbox{for } |x'| \ra 0,      
\end{equation}
and satisfies the condition 
\begin{equation}\label{E2.18} 
\lim_{r\searrow 0} \int_{{\mathcal C}_r} \f T(x)n(x)\g \cdot \vp(x) \,
  dS = 0, 
\end{equation} 
for every $\vp \in C_0^\ue(\R^3,\R^3)$, where ${\mathcal C}_r = \{ (x',
x_3)\in \R^3  \mid |x'| = r \}. $ 
\end{theo}
The proof of this theorem is omitted. \eq{2.15}, \eq{2.16} show that
the stress tensor $T$ is infinitely differentiable on $\R^3\setminus
\ell = \{ (x',x_3)\in \R^3 \mid x'\neq 0\}$. In particular, it is
infinitely differentiable at every point of $\Si$. A simple
computation shows that this is also true for $\na u$. Of course, this
must be the case, since $\Si $ is an artificially introduced surface:
The crystal lattice is undisturbed at this surface.

Let $\Om \subseteq \R^3$ be an open set. To every $C^1$--curve $\ell$
in $\Om$ representing a dislocation curve with unit tangent vector
field $\tau$ and Burgers vector $b$ we define vector and tensor valued
Radon measure $\rho_\ell$, $\hat{b}\otimes \rho_\ell$ respectively, by setting
for all $\vp\in C_0(\Om,\R^3)$, $\tvp \in C_0(\Om,\R^{3\ti3})$,  
\begin{eqnarray}
\label{E2.20}
\lan \rho_\ell,\vp\ran & = & |b|\int_\ell \tau(x)\cdot \vp(x)\,
  ds_x\,,\\ 
\label{E2.21} 
\lan \hb \otimes \rho_\ell,\tvp\ran & = & \int_\ell \f b\otimes
  \tau(x)\g : \tvp(x)\, ds_x\,, 
\end{eqnarray}
with the Nye dislocation tensor $b\otimes \tau(x)$. As usual, for a tensor
valued distribution $w$ and for $\vp\in C_0^\infty(\Om,\R^3)$, $\tvp \in
C_0^\infty(\Om,\R^{3\ti3})$ we define
\[
\lan \div\, w, \vp \ran = - \lan w, \na \vp \ran, \quad 
\lan \ve(w), \vp\ran = \lan w, \ve(\vp) \ran, \quad 
\lan \rot\, w, \tvp \ran = \lan w, \rot\, \tvp \ran. 
\]
\begin{lem}\label{L2.2}
Let $\ell$ be the $x_3$--axis with $\tau(x) = e_3$, and let $u$, $T$
be the functions given in \eq{2.12} -- \eq{2.16} to the Burgers vector
$b = (b_1,0,b_3)$. Then the tensor valued distributions $\tilde{h}_e$
and $T$ defined by 
\begin{eqnarray}
\label{E2.22} 
\lan \tilde{h}_e,\tilde{\vp}\ran & = & \int_{\R^3} \na
  u(x):\tilde{\vp}(x)\, dx,\\
\label{E2.23} 
\lan T, \tilde{\vp} \ran & = & \int_{\R^3} T(x):\tilde{\vp}(x)\, dx
\end{eqnarray}
satisfy the equations 
\begin{eqnarray}
\label{E2.24} -\div\, T & = & 0,\\
\label{E2.25} T & = & D_I\, \ve(\tilde{h}_e),\\
\label{E2.26} \rot\, \tilde{h}_e & = & -\hb \otimes \rho_\ell.
\end{eqnarray}
\end{lem}
Note that the integrals in \eq{2.22}, \eq{2.23} exist because of
\eq{2.17}. $D_I$ is defined in \eq{2.7}. We must omit also the proof
of this lemma.  
\vpn 
We call the two Radon measures $\rho_\ell$ and $\hb\otimes \rho_\ell$
defined in \eq{2.20} and \eq{2.21} vector valued and tensor valued
dislocation measures of the dislocation curve $\ell$. The form of these
measures and the equations \eq{2.22}, \eq{2.26} provide the idea for
the definition of general dislocation measures given now.

Let $\Om \subseteq \R^3$ an open set and let $\mu$ be a scalar Radon measure
on $\Om$. As usual we say that $\mu$ vanishes in a neighborhood ${\cal U}$ of
a point $x$, if $\lan \mu, \vp\ran = 0$ for all $\vp \in C_0({\cal
  U},\R)$. The support ${\rm supp}\,\mu$ of $\mu$ is defined as the set of all
points in $\Om$ which have no neighborhood where $\mu$ vanishes. Let
$\be \in C({\rm supp}\, \mu,\R^3)$ be a given function. Since ${\rm
  supp}\,\mu$ is a relatively closed set in $\Om$, the theorem of
Tietze-Urysohn implies that we can extend $\be$ to a function $\hat{\be}
\in C(\Om,\R^3)$ with $\| \hat{\be}\|_{L^\infty(\Om)} \leq \|
\be\|_{L^\infty({\rm supp}\,\mu)}$. The linear mapping 
\[
\f x \mapsto \vp(x) \g \mapsto \f x\mapsto \hat{\be}(x)\cdot \vp(x) \g:
C_0(\Om,\R^3) \ra C_0(\Om,\R)  
\]
satisfies $\| \hat{\be}\cdot \vp \|_{L^\infty(\Om)} \leq \|
\hat{\be}\|_{L^\infty(\Om)} \| \vp \|_{L^\infty(\Om)}$, hence it is
continuous. Therefore a vector valued Radon measure $\be\mu$ on $\Om$ is
defined by 
\[
\lan \be\mu, \vp \ran = \lan \mu, \hat{\be}\cdot \vp \ran,\qquad \vp \in 
C_0(\Om,\R^3). 
\]
It is not difficult to see that this definition does not depend on the special
continuation $\hat{\be}$ chosen. Therefore the notation $\be \mu$ is
justified. 
\begin{tdefi}\label{D2.3}
The set of all measures $\be \mu$ on $\Om$ with a scalar Radon
measure $\mu$ and $\be \in C({\rm supp}\,\mu,\R^3)$ is denoted by
$\EM(\Om)$. We call $\rho = \tau \mu \in \EM(\Om)$ a vector valued
dislocation measure, if $\mu$ is a nonnegative measure, if $\tau \in
C({\rm supp}\,\mu,\R^3)$ satisfies $|\tau(x)| = 1$ for all $x \in
{\rm supp}\,\mu$ and if a function $h \in L^{1,{\rm loc}}(\Om,\R^3)$
exists such that $\rho = \rot\, h$. The set of vector valued
dislocation measures is denoted by $\EM_d(\Om)$. For $\rho = \tau
\mu \in \EM_d(\Om)$ and $\tvp \in C_0(\Om,\R^{3\ti3})$ we also use
the notation and definition
\[ 
|\rho| = \mu, \quad \frac{\rho}{|\rho|} = \tau,
\quad
\lan \hb \otimes \rho, \tvp \ran = \lan |\rho|, (\hb \otimes \tau) :
\tvp \ran. 
\]
$\hb \otimes \rho$ is called tensor valued dislocation measure. If a 
dislocation measure has a density in $L^1(\Om)$, we call it dislocation
density and denote the measure and the density by the same symbol.
\end{tdefi} 
{\bf Remark} The dislocation measure of a dislocation curve $\ell$
defined in \eq{2.20} belongs to $\EM_d(\Om)$ and has the form $
\rho_\ell = \tau |\rho| $ with $|\rho| = |b|\, {\cal H}_\ell$, where
${\cal H}_\ell = {\cal H}^1 \lfloor \ell$ is the one-dimensional
Hausdorff measure restricted to the curve $\ell$.   

This definition suggests to generalize the problem \eq{2.24} --
\eq{2.26} to a boundary value problem in a domain $\Om \subseteq \R^3$
with an arbitrarily given tensor valued dislocation measure of the form
$-\hb \otimes \rho$ on the right hand side of \eq{2.26}. Before we
formulate this general problem we discuss the meaning of 
$\tilde{h}_e$ in the context of the theory of viscoplasticity at small
strains. In this theory one uses the additive decomposition  
\[
\na u = (\na u - \tilde{h}_p) + \tilde{h}_p
\] 
of the deformation gradient $\na u$ into a plastic part $\tilde{h}_p$
and an elastic part $\na u - \tilde{h}_p$, where only the elastic part
generates the stress field:  
\[
T = D \ve(\na u - \tilde{h}_p).
\]
Comparing this equation with \eq{2.25} we see that $\tilde{h}_e$ is
the elastic part of the deformation gradient:
\begin{equation}\label{E2.29}
\tilde{h}_e = \na u - \tilde{h}_p\,.
\end{equation}
In the following we work with $\tilde{h}_p$ instead of $\tilde{h}_e$.
We eliminate $\tilde{h}_e$ in \eq{2.26} by using \eq{2.29}. If we also
replace $\rho_\ell$ on the right hand side of \eq{2.26} by an
arbitrary dislocation measure $\rho \in \EM_d(\Om)$ we obtain
\begin{equation}\label{E2.30}
\rot\, \tilde{h}_p = - \rot\, \tilde{h}_e = \hb \otimes \rho.
\end{equation}
This equation can be simplified slightly by noting that if the
function $h_p \in L^{1,{\rm loc}}(\Om,\R^3)$ satisfies 
\begin{equation}\label{E2.31}
\rot\, h_p = \rho,
\end{equation}
then \eq{2.30} is fulfilled with $\tilde{h}_p = \hb \otimes
h_p$\,. Taking this expression for $\tilde{h}_p$, inserting \eq{2.29}
into \eq{2.25} and replacing \eq{2.26} by \eq{2.31}, we obtain the
boundary value problem for the displacement field $u$ and the stress
field $T$ in an open set $\Om \subseteq \R^3$ representing the material 
points of a viscoplastic body:   
\begin{eqnarray}
-\div\, T & = & 0, \label{E2.32} \\
T & = & D\big( \ve(\na u) - \ve(\hb \otimes h_p)\big),
  \label{E2.33} \\ 
\rot\, h_p & = & \rho, \label{E2.34} \\
T\rain{\pa\Om} n_B &=& \gm. \label{E2.35}
\end{eqnarray}
The elasticity tensor $D:\ES^3 \ra \ES^3$ can be any linear,
symmetric, positive definite mapping, The dislocation measure $\rho
\in \EM_d(\Om)$ and the boundary data $\gm$ are given, $\hb$ is the
unit vector in direction of the Burgers vector.

We note that the splitting \eq{2.29} of $\tilde{h}_e$ into a gradient
field and a field $\tilde{h}_p$ satisfying $\rot\, \tilde{h}_p = -
\rot\, \tilde{h}_e$ is not unique, since we can add the same gradient field
to $\na u$ and $\tilde{h}_p$ and obtain a new splitting. This means in
particular, that if $(u,T,h_p)$ is a solution of \eq{2.32} --
\eq{2.35} and if $\Gm \in L^{1,{\rm loc}}(\Om,\R)$, then we obtain
another solution $(u',T',h'_p)$ by setting 
\[
u' = u+\hb \Gm,\quad T' = T,\quad h'_p = h_p + \na \Gm,
\]
since these equations imply $\rot\, h'_p = \rot\, h_p = \rho$ and 
\[
\ve(\na u' - \hb \otimes h'_p) = \ve(\na u + \hb \otimes \na \Gm - \hb
\otimes h_p - \hb \otimes \na \Gm) = \ve(\na u - \hb \otimes h_p). 
\]
\section{The evolution equation for dislocation measures}\label{S3} 
Let $[0,T_e)$ be a time interval with $T_e > 0$. In the following we 
formulate 
an evolution equation for dislocation measures
\[
\rho: [0,T_e)\ra \EM_d(\Om)
\]
depending on the time. We base this formulation on four principles:
\begin{itemize}
\item[(P1)] By \refd{2.3}, the dislocation measure $\rho(t)$ and the time
  derivative $\pa_t \rho(t)$ must be a rotation field. Therefore the evolution
  equation must be of the form
\begin{equation}\label{E3.1}
\pa_t \rho = \rot_x \al[T,\rho,b],
\end{equation}
with a function $\al:C(\ov{\Om},\ES^3)\ti\EM_d(\Om)\ti \R^3 \ra \EM(\Om)$ to
be determined. 
\item[(P2)] There must exist a free energy $\psi(\ve,h_p)$ and a flux
$q(\ve,u_t,h_p)$ of the free energy such that the Clausius-Duhem inequality
holds:
\begin{equation}\label{E3.2}
\pa_t  \psi + \div_x  q \leq 0.
\end{equation}
\item[(P3)] The evolution equation \eq{3.1} must allow for solutions $t \ra
  \rho_{\ell(t)}$, which are dislocation measures of dislocation curves, which
  move with driving force given by the Peach-Koehler force  
\begin{equation}\label{E3.3}
F = \tau \ti T b.
\end{equation}
\item[(P4)] Plastic deformation is volume conserving. This means that we must 
  have 
\begin{equation}\label{E3.4}
\mathrm{trace}\, \tilde{h}_p= \mathrm{trace} (\hb \otimes h_p) = 0.
\end{equation}
\end{itemize}
\subsection{Conditions (P1) and (P2)}\label{S3.1}
We first discuss the consequences of (P1) and (P2). Combination of \eq{2.32}
-- \eq{2.35} with \eq{3.1} yields the closed system of partial differential
equations governing the evolution of the dislocation measure and the stress
field in a viscoplastic body:
\begin{eqnarray}
\label{E3.5} -\div_x  T(x,t)  &= & 0, \\
\label{E3.6} T(x,t) & = & D\Ba \ve \f \na_x u(x,t)\g-\ve\f \hb \otimes
     h_p (x,t)\g\Bz,\\ 
\label{E3.7} \rot_x h_p(x,t) &=& \rho(x,t), \\
\label{E3.8} \pa_t \rho(x,t) & = & \rot_x \al[T,\rho,b](x,t), \\ 
\label{E3.9} T(x,t)n_B(x,t) &=& \gm(x,t), \qquad (x,t) \in \pa\Om\ti [0,T_e), 
\end{eqnarray}
where the first four equations must hold for $(x,t) \in \Om\ti [0,T_e)$.  
\begin{lem}\label{L3.1}
(i) Let $\Om \subseteq \R^3$ be an open, bounded, simply connected
set. Then $(u,T,h_p,\rho)$ satisfies the equations \eq{3.5} -- \eq{3.8}
if and only if there is a function $\Gm:\Om \ti [0,T_e) \ra \R$ such
that $(u,T,h_p,\Gm)$ solves the equations 
\begin{eqnarray}
\label{E3.10} -\div_x T  &= & 0, \\
\label{E3.11} T & = & D \f \ve ( \na_x u )-\ve ( \hb\otimes
   h_p) \g,\\ 
\label{E3.12} \pa_t h_p &=& \al[T,\rot_x h_p,b] + \na_x \Gm. 
\end{eqnarray}
On the other hand, if $(u,T,h_p,\Gm)$ solves \eq{3.10} -- \eq{3.12}
  define $\rho = \rot_x h_p$. Then $(u,T,h_p,\rho)$ satisfies \eq{3.5}
  -- \eq{3.8}.  
\vpn
(ii) Let the free energy and the flux be given by 
\begin{eqnarray}
\label{E3.13}
\psi\f \ve(\na_x u),h_p\g &=& \ha \f D \ve(\na_x u - \hb\otimes
  h_p) \g : \ve(\na_x u -\hb\otimes h_p), \\
\label{E3.14}
q\f \ve,u_t.h_p,\Gm\g &=& - T(u_t - \hb \Gm). 
\end{eqnarray}
Then the Clausius-Duhem inequality \eq{3.2} holds for every solution
$(u,T,h_p,\Gm)$ of \eq{3.10} -- \eq{3.12} if $\al$ satisfies for all
points $(T,\rho) \in \ES^3 \ti \R^3$ the inequality 
\begin{equation}\label{E3.15}
(T\hb) \cdot \al[T,\rho,b]\geq 0.
\end{equation}
\end{lem}
{\bf Proof:} Let $(u,T,h_p,\rho)$ be a solution of \eq{3.5} --
\eq{3.8}. Combination of \eq{3.7} and \eq{3.8} yields 
\[
\rot_x \big(  \pa_t h_p - \al \big) = \pa_t \rho - \rot_x \al = 0.   
\]
Since $\Om$ is simply connected, this equation implies that $\pa_t h_p
- \al$ is a gradient field. Consequently there is a function $\Gm:
\Om \ti [0,T_e) \ra \R$ such that $\pa_t h_p - \al[T,\rho,b] =
\na_x \Gm$ holds. From this equation we obtain \eq{3.12} if we use
\eq{3.7} to eliminate $\rho$ in the argument of $\al$. On the other
hand, if $(u,T,h_p,\Gm)$ solves \eq{3.10} -- \eq{3.12} then we obtain
from \eq{3.12} for $\rho = \rot_x h_p$ that
\[
\pa_t \rho = \rot_x \pa_t h_p = \rot_x \al[T,\rho,b],
\]
since $\rot_x \na_x \Gm = 0$. This proves (i). To prove (ii) we infer
from \eq{3.13} and \eq{3.11} that  
\begin{align*}
\pa_t \psi\f \ve(\na_x u),h_p) &= \na_\ve \psi: 
  \ve(\na_x u_t - \hb\otimes \pa_t h_p) = T:\ve(\na_x
  u_t - \hb\otimes \pa_t h_p) \\
& = T:(\na_x u_t - \hb\otimes \pa_t h_p) = T:\na_x u_t -
  (T\hb)\cdot (\pa_t h_p), 
\end{align*}
where we used several times that $T(x,t)$ is a symmetric matrix.
\eq{3.14} yields
\begin{align*}
&\div_x\, q = - \div_x \f T(u_t - \hb \Gm)\g \\
&= -(\div_x\, T^T)\cdot(u_t - \hb\Gm) - T^T : \na_x u_t + T^T:(\hb\otimes
  \na_x\Gm) = -T:\na_x u_t + (T\hb)\cdot \na_x\Gm,
\end{align*}
where we employed that $\div_x T^T = \div_x T = 0$, by
\eq{3.10}. Combination of the last two equations with \eq{3.12}
and \eq{3.15} results in   
\[
\frac{\pa}{\pa t}\, \psi + \div_x q = - (T\hb)\cdot (\pa_t h_p -
  \na_x\Gm) = - (T\hb)\cdot \al \leq 0.
\]
\qed 
\\ 
{\bf Remark} This lemma shows that we are free to choose any function $\Gm$ in
the evolution equation \eq{3.12}. This freedom is used in \cite{ABR08} to
introduce an additional field variable to include dislocation nucleation. In
our investigation we choose for simplicity $\Gm = 0$, which avoids the unusual
term $(T\hb) \Gm$ in the free energy flux \eq{3.14}.
\subsection{Condition (P3)}\label{S3.2}
Next we construct a function $\al$, for which condition (P3) is satisfied. Let
$\Sph^2$ denote the unit sphere in $\R^3$ and let $\tal \in C(\Sph^2 \ti
\R^3,\R^3)$ be a given function. For $T \in C(\ov{\Om},\ES^3)$, $\rho = \tau
|\rho| \in \EM_d(\Om)$ and $b \in \R^3$ we set 
\begin{equation}\label{E3.16}
\al[T,\rho,,b] = \tal (\tau,\tau \ti Tb) \ti \tau\, |\rho|.
\end{equation}
This defines a function $\al:C(\ov{\Om},\ES^3)\ti\EM_d(\Om)\ti \R^3 \ra
\EM(\Om)$, since by \refd{2.3} we have $\tau =
\frac{\rho}{|\rho|} \in C({\rm supp}\, |\rho|,\R^3)$, whence $\tal (\tau,\tau
\ti Tb) \ti \tau \in C({\rm supp}\, |\rho|,\R^3)$, which implies that the
right hand side of \eq{3.16} is in $\EM(\Om)$, again by \refd{2.3}.    
With the notation introduced in \refd{2.3} we write the right hand
side of \eq{3.16} as $\tal \f \frac{\rho}{|\rho|},\frac{\rho}{|\rho|}
\ti Tb \g \ti \rho$.
\begin{tdefi}\label{D3.2}
Let $\al$ be defined by \eq{3.16}, let $b \in \R^3$, and let $T:\Om \ti
[0,T_e) \ra \ES^3$ with $T(t) \in C(\ov{\Om},\ES^3)$ be given. The
time dependent dislocation measure $\rho:[0,T_e) \ra \EM_d(\Om)$ is a
solution of  
\begin{equation}\label{E3.17}
\pa_t \rho = \rot_x \Big(\tal \f \frac{\rho}{|\rho|},\frac{\rho}{|\rho|}
  \ti Tb \g \ti \rho \Big),   
\end{equation}
if for all $\vp \in C^\infty_0 \big((0,T_e) \ti \Om,\R^3\big)$ the integrals
in the following equation exist and satisfy
\[
- \int_0^{T_e} \big\lan \rho(t), \pa_t \vp(t) \big\ran\, dt = \int_0^{T_e}
\Big\lan \tal \f \frac{\rho}{|\rho|},\frac{\rho}{|\rho|}
  \ti Tb \g \ti \rho, \rot_x \vp(t) \Big\ran\, dt.
\]
\end{tdefi}
As a consequence of the next theorem we see that condition (P3) is
satisfied for the evolution equation \eq{3.17}. We need two
definitions to state this theorem. Let the family $t \mapsto \ell(t)$ represent
a moving dislocation curve with tangent vector $\tau(x,t)$ at the point
$x \in \ell(t)$. By ${\rm proj}_{\tau(x,t)}$ we denote the orthogonal
projection
\begin{equation}\label{E3.18}
{\rm proj}_{\tau(x,t)} : \R^3 \to H(x,t) = \{ \xi \in \R^3 \mid \xi
\cdot \tau(x,t)=0\} \subseteq \R^3 
\end{equation}
to the orthogonal space of $\tau(x,t)$. To define the normal velocity
$v(x_0,t_0)$ of the dislocation curve at time $t_0$ at $x_0 \in
\ell(t_0)$, let $x(t)$ be the intersection point of $\ell(t)$ with
$H(x_0,t_0)$. Set
\[
v(x_0,t_0) = \frac{d}{dt} x(t)\rain{t=t_0}. 
\]
\begin{theo}\label{T3.3}
Let $\al$ be defined by \eq{3.16} with a given function $\tal \in C(\Sph^2 \ti 
\R^3,\R^3)$. \\ 
(i) $\al$ satisfies the inequality \eq{3.15} if   
\begin{equation}\label{E3.19}
\xi \cdot \tal (\tau,\xi) \geq 0, \quad \mbox{for all } (\tau,\xi) \in \Sph^2
\ti \R^3. 
\end{equation}
(ii) Assume that $T \in C \f [0,T_e] \ti \ov{\Om},\ES^3 \g$ is a given
stress field and that $\rho_{\ell(t)} = \tau\, |\rho_{\ell(t)}| \in
\EM_d(\Om)$ is the dislocation measure of a dislocation curve
$\ell(t)$ to the Burgers vector $b\in \R^3$. For $x \in \ell(t)$ let
$F(x,t) = \tau(x,t) \ti T(x,t)b$ be the Peach-Koehler force. Then
$\rho_{\ell(t)}$ solves the evolution equation \eq{3.17}, if and only
if for every $x \in \ell(t)$ the normal speed is
\begin{equation}\label{E3.20}
v(x,t) =  {\rm proj}_{\tau(x,t)} \tal \Ba \tau(x,t),F(x,t)\Bz.  
\end{equation}
\end{theo}
We must omit the proof of this theorem. \eq{3.20} shows that $\tal$
can be considered to be a constitutive function determining the
relation between the normal velocity and the Peach-Koehler force $F$.
Since the continuity of $\tal$ and the inequality \eq{3.19} imply
$\tal(\tau,0) = 0$, the normal velocity is equal to zero if $F=0$.
Therefore we can consider $F$ to be the driving force for the movement
of $\ell(t)$. In this sense, condition (P3) is satisfied by the
evolution equation \eq{3.17}. The evolution equation \eq{3.17} was in
principal derived in \cite{Mura63}, following ideas in \cite{Kroe58}, cf. also 
\cite{Fox66}.   

With the function $\al$ determined in \eq{3.16}, with $\Gm=0$ and with
$\rot_x h_p = \tau\, |\rot_x h_p| \in \EM_d(\Om)$ the system \eq{3.10}
-- \eq{3.12} combined with the boundary condition \eq{3.9} takes the
form
\begin{eqnarray}
\label{E3.21}  -\div_x T &=& 0,  \\
\label{E3.22} T &=&  D \f \ve ( \na_x\, u) - \ve( \hb \otimes h_p ) \g,  
  \\  
\label{E3.23} \pa_t h_p &=& \tal \f \tau, \tau \ti Tb \g \ti \rot_x
  h_p \,,
  \\ 
\label{E3.24} T\rain{\pa \Om} n_B &=& \gm,\qquad \mbox{on } \pa\Om \ti
   [0,T_e).   
\end{eqnarray}
\subsection{Condition (P4)}\label{S3.3}
To determine the consequences of condition (P4) we differentiate
\eq{3.4} with respect to $t$ and use \eq{3.23} to compute  
\[
0 = {\rm trace} (\hb\otimes \pa_t h_p) = \hb \cdot \pa_t h_p = \hb
\cdot (\tal \ti \tau) |\rot h_p| = 0. 
\] 
Here we used that $\rot_x h_p = \tau |\rot h_p|$. From this equation
we infer that  
\begin{equation}\label{E3.25}
b \cdot \Ba \tal \f \tau(x,t),\tau(x,t) \ti T(x,t)b \g \ti \tau(x,t)
  \Bz = 0.  
\end{equation}
must hold for all $(x,t) \in \Om \ti [0,T_e)$ if condition (P4) is
satisfied. This equation has the following consequence for the
movement of dislocation curves: 
\begin{coro}\label{C3.4}
Let the stress field $T:\Om \ti [0,T_e) \ra \ES^3$ with $T(t) \in
C(\ov{\Om},\ES^3)$ be given, and let $\rho_{\ell(t)} = \tau\,
|\rho_{\ell(t)}|$ be the dislocation measure of a moving dislocation
curve $\ell(t)$. If $\rho_{\ell(t)}$ solves the evolution equation
\eq{3.17} and if the condition \eq{3.25} holds for a point $x \in
\ell(t)$ at which $\tau(x,t)$ is not parallel to $b$, then the
normal speed $v(x,t)$ of the dislocation curve is parallel to the
vector ${\rm proj}_{\tau(x,t)} b$.
\end{coro}
{\bf Proof:} $\tau$ and $b$ are parallel if and only if $b \ti \tau = 0$. If
this product is not zero, then \eq{3.25} implies that the vector
$\tal(\tau,\tau \ti Tb)$ must lie in the plane spanned by $\tau$ and $b$. In
this case the vectors ${\rm proj}_{\tau(x,t)} b$ and ${\rm proj}_{\tau(x,t)}
\tal(\tau,\tau \ti Tb)$ lie on the same line, so they are parallel. Since by
\reft{3.3}(ii) we have $v = {\rm proj}_{\tau(x,t)} \tal(\tau,F)$, the
corollary follows.  \qed \vpn {\bf Remark} If the tangent vector $\tau(x,t)$
is parallel to the Burgers vector $b$, then $\ell(t)$ is called a screw
dislocation at $x \in \ell(t)$. The condition that the volume is not changed
by plastic deformation thus requires that if $\ell(t)$ is not a screw
dislocation at $x \in \ell(t)$, then the normal speed $v(x,t)$ must be a
linear combination of $b$ and $\tau(x,t)$. Only for screw dislocations the
direction of the normal speed is not determined by the Burgers
vector. However, a dislocation curve can move freely in any direction only if
it is a screw dislocation at every point, which means that it must be a
straight dislocation line in direction of the Burgers vector, and if the
movement is a parallel shift. As soon as the movement ceases to be a parallel
shift, points appear on $\ell(t)$, at which $\tau$ and $b$ are not parallel,
restricting the direction of the movement of $\ell(t)$ at these points.

From \eq{3.25} we see that in order to guarantee that condition (P4)
is satisfied, the function $\tal$ must be such that
\[
b \cdot \f \tal(\tau,\xi) \ti \tau \g = 0,
\]
for all $(\tau,\xi) \in \Sph^2 \ti \R^3$. If $b$ and $\tau$ are
linearly independent, this holds if and only if $\tal(\tau,\xi)$
belongs to the linear span of $b$ and $\tau$, hence there are real
valued functions $f_1$ and $f_2$ such that $\tal(\tau,\xi) = b
f_1(\tau,\xi) + \tau f_2(\tau,\xi)$. This implies $\tal(\tau, \xi )
\ti \tau = (b \ti \tau) f_1(\tau,\xi)$, hence $\tau f_2$ does not
contribute to the evolution equation and we can omit this term and
construct $\tal$ in the form
\[
\tal(\tau,\xi) = b f_1(\tau,\xi),
\]
where the function $f_1$ must be chosen such that the inequality
\eq{3.19} holds. Under many possibilities a simple choice is  
\begin{equation}\label{E3.26}
\tal(\tau,\xi) = \frac{b}{|b \ti \tau|}\, f\Ba \frac{b}{|b \ti \tau|}
  \cdot   \xi \Bz, 
\end{equation}
with a function $f \in C(\R,\R)$ satisfying $r \cdot f(r) \geq 0$ for
all $r \in \R$. We insert \eq{3.26} into \eq{3.17} and obtain
the evolution equation
\begin{equation}\label{E3.27}
\pa_t \rho = \rot_x \left( f\Ba \frac{b \ti
\tau}{|b \ti \tau|} \cdot Tb \Bz \frac{b \ti \tau}{|b \ti \tau|}
  |\rho| \right), 
\end{equation}
for $\rho = \tau |\rho| \in \EM_d(\Om)$. Here we used that $b \cdot (\tau \ti
Tb) = (b \ti \tau) \cdot Tb$ to rewrite the argument of $f$. The choice of the
$\tal$ in \eq{3.26} is justified by the following result:
\begin{theo}\label{T3.5}
(i) For the function $\tal$ defined in \eq{3.26} the inequality
\eq{3.19} holds. \\
(ii) 
Assume that $T \in C \f [0,T_e] \ti \ov{\Om},\ES^3 \g$ is a given stress
field and that $\rho_{\ell(t)} = \tau\, |\rho_{\ell(t)}| \in \EM_d(\Om)$ is
the dislocation measure of a dislocation curve $\ell(t)$ to the Burgers vector 
$b\in \R^3$. Let $F = \tau \ti Tb$ be the Peach-Koehler force, and for $x \in
\ell(t)$ denote the unit vector in the direction of ${\rm proj}_{\tau(x,t)}$
by   
\[
b^{\bot\tau}(x,t) = \frac{{\rm proj}_{\tau(x,t)} b}{|{\rm
    proj}_{\tau(x,t)} b|}. 
\]
Then $\rho_{\ell(t)}$ solves the evolution equation \eq{3.27}, if and
only if $\ell(t)$ moves with the normal speed
\begin{equation}\label{E3.28}
v(x,t) =  f \f b^{\bot\tau}(x,t) \cdot F(x,t) \g\, b^{\bot\tau}(x,t). 
\end{equation}
(iii) Let $(u,T,h_p)$ be a solution of the system 
\begin{eqnarray}
\label{E3.29}  -\div_x\, T &=& 0,  \\
\label{E3.30} T &=&  D \f \ve ( \na_x\, u) - \ve( \hb \otimes h_p ) \g,  
\\  
\label{E3.31} \pa_t h_p &=& f\Ba \frac{b \ti \tau}{|b \ti \tau|} \cdot
      Tb \Bz \frac{b \ti \tau}{|b \ti \tau|} |\rot_x h_p| 
\end{eqnarray}
in $\Om \ti [0,T_e)$, where $\rot_x h_p = \tau |\rot_x h_p|$. If at
time $t=0$ we have $b \cdot h_p(x,0) = 0$ for all $x \in \Om$, then
condition \eq{3.4} holds for all $(x,t) \in \Om \ti [0,T_e)$.  
\end{theo}
{\bf Remarks} \eq{3.28} shows that under the choice of $\tal$
given in \eq{3.26} the driving force for the movement of dislocation
curves is the component $b^{\bot\tau} \cdot F$ of the Peach-Koehler
force in the direction of the vector ${\rm proj}_\tau b$. This justifies
\eq{3.26}. The equation \eq{3.31} is obtained by insertion of the function
$\tal$ from \eq{3.26} into \eq{3.23}, using again that $b \cdot (\tau \ti Tb) 
= (b \ti \tau) \cdot Tb$.  
\\[1ex] 
{\bf Proof:} For the proof of (i) note that the assumption $r \cdot
f(r) \geq 0$ implies 
\[
\xi \cdot \tal(\tau, \xi ) = \frac{\xi \cdot b}{|b \ti \tau|}\,
f \Ba \frac{b \cdot \xi}{|b \ti \tau|} \Bz \geq 0, \quad (\tau,\xi) \in
\Sph^2 \ti \R^3.  
\]
To verify (ii) note that since the evolution equation \eq{3.27} is
obtained from \eq{3.17} by insertion of the function $\tal$ defined in
\eq{3.26}, we obtain from \reft{3.3}(ii) that $\rho_{\ell(t)}$ is a
solution of \eq{3.27} if and only if the dislocation curve moves with
normal velocity
\begin{equation}\label{E3.32}
v = {\rm proj}_{\tau} \left( \frac{b}{|b \ti \tau|}\, f\Ba \frac{b}{|b
    \ti \tau|} \cdot F \Bz \right) = f\Ba \frac{b}{|b
    \ti \tau|} \cdot F \Bz \frac{ {\rm proj}_{\tau} b}{|b \ti \tau|}.     
\end{equation}
Since $|b \ti \tau| = |{\rm proj}_\tau b|$, we have $\frac{ {\rm
proj}_{\tau} b}{|b \ti \tau|} = b^{\bot\tau}$. Noting that $F = \tau \ti Tb$
is orthogonal to $\tau$, we thus conclude  
\[
\frac{1}{|b \ti \tau|}\, b \cdot F = \frac{1}{|b \ti \tau|}
  ({\rm proj}_{\tau} b ) \cdot F = b^{\bot\tau} \cdot F\,.     
\] 
\eq{3.28} is obtained by insertion of this equation into \eq{3.32}. 

To prove (iii) observe that the right hand side of \eq{3.31} is
orthogonal to the vector $b$, which implies that $\pa_t (\hb \cdot
h_p) = \hb \cdot \pa_t h_p = 0$. This equation and the assumption $b
\cdot h_p(x,0) = 0$ together imply $\hb \cdot h_p(x,t) = 0$   
for all $(x,t) \in \Om \ti [0,T_e)$, from which we conclude that ${\rm
  trace}(\hb \otimes  h_p) = \hb \cdot h_p = 0$. This is \eq{3.4}. \qed 
\vpn
It is possible to simplify \eq{3.31} slightly by working with a suitable
coordinate system. To show this we assume that $T$ is a given stress
field and that $h_p$ is a solution of \eq{3.31} satisfying $\hb \cdot
h_p = {\rm trace}(\hb \otimes h_p) = 0$. Choose the cartesian
coordinates $(x_1,x_2,x_3)$ such that the $x_2$--axis points into the
direction of the Burgers vector $b$ and let $h_p =
(h_{p1},h_{p2},h_{p3})$ be the components of $h_p$ in this coordinate 
system. Since $\hb = (0,1,0)$, the equation $\hb \cdot h_p = 0$ is
equivalent to $h_{p2} = 0$. Using this property of $h_p$ and noting that $\tau
= \frac{\rot_x h_p}{|\rot_x h_p|}$, we obtain by a computation that in these
coordinates 
\begin{align*}
&\rot_x h_p = \f \pa_{x_2} h_{p3}, \pa_{x_3} h_{p1} - \pa_{x_1} h_{p3},
  -\pa_{x_2} h_{p1} \g, \qquad \hb \ti \rot_x h_p = - \pa_{x_2} h_p.
\\
&\frac{b \ti \tau}{|b \ti \tau|} = \frac{b \ti  \rot_x h_p } 
   {|b \ti  \rot_x h_p  |} = - \frac{ \pa_{x_2} h_p} {|\pa_{x_2} h_p|},  
  \qquad
    |\rot_x h_p| = \sqrt{|\pa_{x_2} h_p|^2 + (\rot_2 h_p)^2 },
\end{align*}
where $\rot_2 h_p = \pa_{x_3} h_{p1} - \pa_{x_1} h_{p3}$\,. With these
expressions \eq{3.31} takes the form
\begin{equation}\label{E3.33}
\pa_t h_p = - f \Ba - \frac{ \pa_{x_2} h_p} {|\pa_{x_2} h_p|} \cdot Tb
  \Bz \frac{ \pa_{x_2} h_p} {|\pa_{x_2} h_p|} \sqrt{|\pa_{x_2} h_p|^2 +
    (\rot_2 h_p)^2 }.  
\end{equation}
This is a system of two equations, since $h_{p2}$ vanishes
identically.   
\section{Dislocations moving in slip planes}\label{S4}
We finally consider the situation, where every dislocation curve is
contained in a plane and moves within this plane. Different
dislocation curves can be contained in the same plane or in parallel
planes. The planes are called slip planes. 

Let $g \in \R^3$ be a unit vector, which is normal to all slip planes.
Every tangent vector $\tau$ to a dislocation curve in a slip plane is
normal to the vector $g$. Since the dislocation curve moves in the
slip plane, the normal velocity $v$ of the dislocation curve must also
be normal to $g$. By \reft{3.5}(ii), the vector $v$ is parallel to the
vector ${\rm proj}_\tau b$, which is the orthogonal projection of the
Burgers vector $b$ to the orthogonal space of $\tau$. This implies
that $b$ itself is orthogonal to $g$. Therefore the vectors $\tau$ and
$b$ span the two dimensional subspace parallel to the slip planes. The
vector $b \ti \tau$ is normal to this subspace, hence
\[
\frac{b \ti \tau}{|b \ti \tau|} = \pm g.
\]
We insert this equation into \eq{3.27} and obtain the evolution equation
\begin{equation}\label{E4.1} 
\pa_t \rho = \rot_x \Ba {\pm} g\, f \f {\pm} g \cdot Tb \g |\rho| \Bz.
\end{equation}
The change of sign can occur at points, where the dislocation curve is
a screw dislocation, that is at points where $\tau$ and $b$ are
parallel. We have not accounted for this situation in the definition
\eq{3.26} of the constitutive function $\tal$, which is used to obtain
\eq{3.31}. From \refc{3.4} we know that at these points the direction
of the normal velocity is not determined by $\tau$ and $b$, and a
separate definition of $\tal$ should be given for such points.
However, under our present assumption that dislocations move in slip
planes, we can obtain an evolution equation valid everywhere by simply
taking the $+$--signs in \eq{4.1}. We avoid to discuss the
justification of this choice in general, but simply assume in the
following that $f$ is chosen as an odd function, in which case the
$\pm$--signs in \eq{4.1} can obviously be replaced by $+$--signs.
\begin{coro}\label{C4.1}
Assume that $T \in C \f [0,T_e] \ti \ov{\Om},\ES^3 \g$ is a given
stress field and that $\rho_{\ell(t)} =
\tau\, |\rho_{\ell(t)}| \in \EM_d(\Om)$ is the dislocation measure of
a dislocation curve $\ell(t)$ to the Burgers vector $b\in \R^3$, which
for all $t \in [0,T_e)$ is contained in a plane normal to the vector
$g$. Then $\rho_{\ell(t)}$ is a solution of the evolution equation
\begin{equation}\label{E4.2}
\pa_t \rho = \rot_x \Ba g\, f \f g \cdot Tb \g |\rho| \Bz ,
\end{equation}
if and only if it moves with the normal speed given in \eq{3.28}. 
\end{coro}
This corollary follows immediately from \reft{3.5}(ii) and the 
construction of \eq{4.2} given above.  

We can simplify \eq{4.2} by introducing a cartesian coordinate system
with the $x_2$--axis pointing into the direction of $b$ and the
$x_3$--axis pointing into the direction of $g$. For a scalar function
$w$ we then have $\rot_x (wg) = -(g \ti \na_x) w = \na_g^\bot w$, with
$\na_g^\bot = (\pa_{x_2}, -\pa_{x_1},0)^T$. Equation \eq{4.2} thus
becomes
\begin{equation}\label{E4.3}
\pa_t \rho = \na_g^\bot \Ba f ( g \cdot Tb ) |\rho| \Bz. 
\end{equation}
The equation for $h_p$ corresponding to the evolution equation
\eq{4.2} is obtained by insertion of the vector $g$ for $\frac{b \ti
  \tau}{|b \ti \tau|}$ in \eq{3.31}. We obtain
\begin{equation}\label{E4.4}
\pa_t h_p = g\, f \f g \cdot Tb \g |\rot_x h_p|. 
\end{equation}
From this equation we see that $\pa_t h_p$ is a scalar multiple of the
vector $g$. We can assume that this is also true for the initial data
$h_p(x,0)$, from which we conclude that there is a function $\ve_p:
\Om \ti [0,T_e) \ra \R$ such that for all $(x,t) \in \Om\ti [0,T_e)$
\begin{equation}\label{E4.5}
h_p (x,t) = \ve_p(x,t)\,g .
\end{equation}
With the coordinate system chosen as above and with the tangential
gradient $\na_g = (\pa_{x_1}, \pa_{x_2},0)^T$ of the slip plane we
obtain from \eq{4.5} that 
\begin{equation}\label{E4.6}
\rot_x h_p = \rot_x (\ve_p g) = \na_g^\bot \ve_p, \qquad
  |\rot_x (\ve_p g)| =  |\na_g^\bot \ve_p | = |\na_g \ve_p|. 
\end{equation} 
We insert the second expression and \eq{4.5} into \eq{4.4} and obtain the
evolution equation for $\ve_p$:
\[
\pa_t \ve_p = f(g \cdot Tb)\, |\na_g \ve_p|.
\]
For dislocations moving in slip planes this equation replaces the third
equation in the system \eq{3.29} -- \eq{3.31}. The model equations for
plastic deformation of materials with dislocations moving in slip
planes thus take the form 
\begin{eqnarray}
\label{E4.7}  -\div_x\, T &=& 0,  \\
\label{E4.8} T &=&  D \f \ve ( \na_x\, u) - m\ve_p \g,  
\\  
\label{E4.9} \pa_t \ve_p &=& f(|b|\, m:T)\, |\na_g \ve_p|,
\\
\label{E4.10} T\rain{\pa \Om} n_B &=& \gm,\qquad \mbox{on } \pa\Om \ti
  [0,T_e), 
\end{eqnarray}
where $m$ is defined in \eq{1.4} and where we used that $\ve( \hb \otimes h_p)
= \ve( \hb \otimes \ve_p g) = m \ve_p$ and $g \cdot Tb = |b|\,m:T$. Of course,
this system has to be supplemented by a boundary condition and an initial
condition for the function $\ve_p$. We do not discuss this here.
\paragraph{Remark} 
Let $(u,T,\ve_p)$ be a solution of \eq{4.7} -- \eq{4.9}. By \eq{3.7}
the dislocation measure to this solution is given by $\rho = \rot_x
h_p \in \EM_d(\Om)$, from which we obtain by \eq{4.6} that $\rho =
\tau |\rho| = \na_g^\bot \ve_p$. Since $\na_g^\bot \ve_p$ is
orthogonal to $g = (0,0,1)^T$, we conclude that also the tangential
vector field $\tau$ of the dislocation measure $\rho$ is orthogonal to
$g$, and this implies that the dislocation curves corresponding to
solutions of the system \eq{4.7} -- \eq{4.9} lie in slip planes $x_3 =
c$ and that the function $(x_1,x_2) \mapsto \ve_p(x_1,x_2,c,t): \R^2
\ra \R$ plays the role of a stream function of the dislocation measure
$\rho$ on the slip plane.

If $\rho = \rho_\ell = \tau |\rho_\ell|$ is the dislocation measure of a
dislocation curve, then $\ve_p$ has a jump along the dislocation curve. By
definition of $\rho_\ell$ in \eq{2.20}, the height of the jump is equal to the
absolute value $|b|$ of the Burgers vector. If $\rho$ is a dislocation
density, then $\ve_p$ is smooth and the level curves of $(x_1,x_2) \mapsto
\ve_p(x_1,x_2,c,t)$ are the averaged dislocation curves.
\vpn
{\bf Summary} Both of the systems \eq{3.29}, \eq{3.30}, \eq{3.33} and \eq{4.7}
-- \eq{4.9} with the corresponding evolution equations \eq{3.27} and \eq{4.3},
respectively, satisfy the conditions (P1) -- (P4). The simpler system \eq{4.7}
-- \eq{4.9} models dislocation curves which lie in slip planes. If this
restrictive assumption is not made, then the more complicated equations
\eq{3.29}, \eq{3.30}, \eq{3.33} must be used.

The system \eq{4.7} -- \eq{4.9} differs from the standard plasticity model
\eq{1.1} -- \eq{1.3} only by the term $|\na_g \ve_p|$. From \eq{4.6} we see
that $|\rho| = |\na_g \ve_p|$, hence $|\na_g \ve_p|$ is the total variation
measure to the dislocation measure and represents the ``absolute value''
of the dislocation density. If this density is high throughout the material,
then $|\na_g \ve_p|$ can be replaced approximately by a constant. In this case
the system \eq{4.7} -- \eq{4.9} reduces to the standard plasticity model
\eq{1.1} -- \eq{1.3}. However, if this assumption is not valid, then \eq{4.7}
-- \eq{4.9} should be used as model for plastic deformation. 


\begin{thebibliography}{99}
\itemsep0pt
{\small 
\bibitem{ABR08} A. Acharya, A. Beaudoin, R. Miller. New perspectives in
  plasticity theory. Dislocation Nucleation, waves, and partial continuity of
  plastic strain rate. {\it Math. Mech. Solids} {\bf 13} (2008), 292--315.  
\bibitem{Acharya11} A. Acharya. Microcanonical entropy and mesoscale
  dislocation mechanics and plasicity. {\it J. Elast.} {\bf 104} (2011),
  23--44  
\bibitem{AlChel07} H.-D. Alber, K. Che\l mi\'nski. Quasi-static
  problems in viscoplasticity theory II: Models with nonlinear
  hardening. {\it Math. Models Meth. Appl. Sci.} {\bf 17,2} (2007), 
  189--213. 
\bibitem{Fox66} N. Fox. A continuum theory of dislocations for single
  crystals. IMA J. Appl. Math. {\bf 2} (1966), 285--298. 
\bibitem{Kos79} A.M. Kosevich. Crystal dislocations and the theory of
  elasticity. {\it In: F.R.N. Nabarro [ed.] Dislocations in Solids }
  33--141, North-Holland 1979.     
\bibitem{Mura63} T. Mura. Continuous distribution of moving dislocations. {\it
  Philosophical Magazine} {\bf 89} (1963), 843--857.
\bibitem{Kroe58} E. Kr\"oner. Kontinuumstheorie der Versetzungen und
  Eigenspannungen. Ergebnisse der Angewandten Mathematik {\bf 5}, Springer
  1958  
} 
\end{thebibliography}
\end{document}